\documentclass{amsart}

\usepackage{amssymb}
\usepackage[all]{xy}
\usepackage{hyperref}

\usepackage{enumitem}   

\setlist[enumerate]{itemsep=.2em,topsep=.2em,leftmargin=1.25em,itemindent=2.0em}

\newtheorem{thm}{Theorem}

\newtheorem{lem}[thm]{Lemma}
\newtheorem{cor}[thm]{Corollary}

\newtheorem{prop}[thm]{Proposition}

\newtheorem{principle}[thm]{Principle} 

\theoremstyle{definition}
\newtheorem{defn}[thm]{Definition}

\newtheorem{say}[thm]{}
\newtheorem{exmp}[thm]{Example}

\newtheorem{rem}[thm]{Remark}          

\newtheorem*{ack}{Acknowledgments}      
\newtheorem{notation}[thm]{Notation}   
  
\newtheorem{defn-thm}[thm]{Definition--Theorem}  
\newtheorem{defn-lem}[thm]{Definition--Lemma}  

\theoremstyle{remark}

\renewcommand{\c}[0]{{\mathbb C}}  

\renewcommand{\o}[0]{{\mathcal O}} 

\renewcommand{\a}[0]{{\mathbb A}}

\newcommand{\p}[0]{{\mathbb P}}

\newcommand{\q}[0]{{\mathbb Q}}
\newcommand{\map}[0]{\dasharrow}
\newcommand{\qtq}[1]{\quad\mbox{#1}\quad}

\newcommand{\im}[0]{\operatorname{im}}    
\newcommand{\proj}[0]{\operatorname{Proj}}

\newcommand{\aut}[0]{\operatorname{Aut}}

\newcommand{\tsum}[0]{\textstyle{\sum}}

\newcommand{\bir}[0]{\operatorname{Bir}}

\def\into{\DOTSB\lhook\joinrel\to}

\def\loccoh#1.#2.#3.#4.{H^{#1}_{#2}(#3,#4)}

\DeclareMathAlphabet{\mathchanc}{OT1}{pzc}%
                                {m}{it}

\newcommand{\gm}[0]{{\mathbb G}_m}
\newcommand{\ga}[0]{{\mathbb G}_a}

\newcommand{\simb}[0]{\stackrel{bir}{\sim}}

\newcommand{\GL}{\mathrm{GL}}
\newcommand{\PGL}{\mathrm{PGL}}
\newcommand{\PSL}{\mathrm{PSL}}

\newcommand{\SL}{\mathrm{SL}}

\usepackage[all]{xy}\xyoption{dvips}

\begin{document}
\bibliographystyle{amsalpha}

\hfill\today

\title{Automorphisms of unstable ${\mathbb P}^1$-bundles}
        \author{J\'anos Koll\'ar}

\begin{abstract}
    Let $P\to X$ be a ${\mathbb P}^1$-bundle over a variety $X$.
The aim of this note is to understand all  connected,  
algebraic groups
$$
\operatorname{Aut}^\circ(P)\subset G\subset \operatorname{Bir}( X\times {\mathbb P}^1).
$$
We get a quite complete answer if $\operatorname{Aut}^\circ(X)$ is a
maximal, connected,  
algebraic subgroup of $\operatorname{Bir}(X)$, and $P$ is sufficiently unstable.

This gives examples of connected,  
algebraic subgroups of $\operatorname{Bir}({\mathbb P}^4)$ that are not contained in any maximal one.
\end{abstract}

   \maketitle

   Enriques described all connected, algebraic  subgroups of the 2-dimensional Cremona group $\bir(\p^2)$ in 1893.
   In dimension 3, Umemura  classified all such subgroups,
   in a series of papers completed in \cite{MR0945907}.
   They are all contained in a maximal one.
Much less is known about algebraic  subgroups of $\bir(\p^n)$ for  $n\geq 4$.

Fong and Zikas  gave examples of connected,   algebraic subgroups of
  $\bir(C\times \p^n)$
that are  not contained in any maximal one, where $C$ is a non-rational curve
 \cite{fong-1, MR4538552}. 
 Fanelli, Floris and Zimmermann constructed similar  subgroups in  
  $\bir(\p^n)$ for  $n\geq 5$, answering a question of Blanc \cite{MR4537333}; see \cite{ffz}.
 The proof in \cite{ffz} relies on the existence of a
 stably rational but nonrational conic bundle 3-fold $X$ \cite{MR86m:14009}, 
 a careful investigation of the  birational geometry of
   unstable $\p^1$-bundles over $X$, and a study of  unipotent group actions.
   For background and previous results, see
   \cite{MR4545346, MR4537333} for birational group actions, and
   \cite{Sar81, MR1311348, km-book, ksc}  for minimal models and birational maps.

   In this note we follow \cite{fong-1, MR4538552, ffz}, but,
   instead of using birational geometry, 
   we focus entirely on the action of the  unipotent radical of an algebraic group. 
   This leads to  examples in $\bir(\p^4)$ as well.

     As another consequence,
our results apply  over  algebraically closed fields of arbitrary  characteristic.

 The main results are Theorems~\ref{main.G.thm} and \ref{best.reg.G.prop.2}.
 Here we state a special case that is easy to formulate,  yet leads to many interesting subgroups of  $\bir(\p^4)$.

        \begin{thm}\label{main.thm} Let $X$ be a smooth, projective variety such  that  $\aut^\circ(X)$ is  a maximal, connected,  
          algebraic subgroup of  $\bir(X)$.
                  Let $L$ be an $\aut^\circ(X)$-linearized, very ample line bundle on $X$, 
          and
          $$
          \aut^\circ\bigl(\p_X(\o_X{\oplus}L)\bigr)\subset  G\subset \bir(X\times \p^1)
          $$
          a connected,  algebraic subgroup.
Assume that   $(X, L)\not\cong \bigl(\p^n, \o_{\p^n}(1)\bigr)$.

          Then there is an  $\aut^\circ(X)$-invariant, effective divisor $ D$   on $X$ such that
          $$
           G\subset \aut^\circ\bigl(\p_X\bigl(\o_X{\oplus}L( D)\bigr)\bigr).
          $$
        \end{thm}

       Note that  \cite{fong-1, MR4538552, ffz} consider cases for which
       $\aut^\circ(X)=\{1\}$, and use a Noether-Fano type argument to show that $G$ preserves the negative section of $\p_X(\o_X{\oplus}L)$.
       It is quite likely that one can prove Theorem~\ref{main.thm}
       using the methods of \cite{fong-1, MR4538552, ffz}, at least in characteristic 0.

        
        For Corollaries~\ref{main.cor.2}--\ref{main.cor.1}
        we use the notation and assumptions of Theorem~\ref{main.thm}, and the description of $\aut^\circ\bigl(\p_X(\o_X{\oplus}L)\bigr)$ given in
        Paragraph~\ref{aut.pb.say}.

       \begin{cor}\label{main.cor.2}
$\aut^\circ\bigl(\p_X(\o_X{\oplus}L)\bigr)$ is a  maximal, connected, algebraic subgroup of
$\bir(X\times \p^1)$ iff
$H^0(X, L)=H^0\bigl(X, L(D)\bigr)$ 
for every   $\aut^\circ(X)$-invariant  divisor $D$ on $X$. \qed
       \end{cor}

       \begin{cor}\label{main.cor.2.a}
$\aut^\circ\bigl(\p_X(\o_X{\oplus}L)\bigr)$ is a  maximal, connected, algebraic subgroup of
$\bir(X\times \p^1)$ for every very ample $L$ iff  there are no
       $\aut^\circ(X)$-invariant  divisors on $X$. 
       \end{cor}

       Proof. One direction follows from Corollary~\ref{main.cor.2}.
       Conversely, if $D\neq 0$ is an $\aut^\circ(X)$-invariant  divisor
       then there is a very ample $L$ such that
       $h^0(X, L)<h^0\bigl(X, L(D)\bigr)$.
       Then $\aut^\circ\bigl(\p_X(\o_X{\oplus}L)\bigr)\subsetneq \aut^\circ\bigl(\p_X(\o_X{\oplus}L(D))\bigr)$\qed 
       \medskip

        If there are few  $\aut^\circ(X)$-orbits, quite likely this can  also be  proved using the Noether-Fano method \cite{ksc}, but the general case is not clear to me.
        
        \begin{cor}\label{main.cor.1}
          $\aut^\circ\bigl(\p_X(\o_X{\oplus}L)\bigr)$ is not contained in any maximal, connected, algebraic subgroup of
  $\bir(X\times \p^1)$, if 
                      there is an $\aut^\circ(X)$-invariant divisor $D'$ such that
          $h^0\bigl(X, \o_X(D')\bigr)>1$.
              \end{cor}

       Proof.  By Theorem~\ref{main.thm},
if $\aut^\circ\bigl(\p_X(\o_X{\oplus}L)\bigr)\subset  G$, then 
$ G\subset \aut^\circ\bigl(\p_X(\o_X{\oplus}L( D))\bigr)$ for some $\aut^\circ(X)$-invariant, effective divisor $ D$.

        Since  $h^0\bigl(X, \o_X(D')\bigr)>1$, we see that 
        $h^0\bigl(X, L(D+D')\bigr)>h^0(X, L( D))$. 
        Thus   $\aut^\circ\bigl(\p_X(\o_X{\oplus}L( D))\bigr)$ is contained in the  strictly bigger group
        $\aut^\circ\bigl(\p_X(\o_X{\oplus}L(D+D'))\bigr)$. \qed

        \medskip

        Using Umemura's classification of maximal, connected, algebraic subgroups of $\bir(\c\p^3)$,  there are 3 cases  over $\c$ to which Corollary~\ref{main.cor.1} applies.
\cite{MR4537333} gives an algebraic proof of  the  classification, which  also works in characteristics $p>5$ (though some parts are not yet fully written down).

\begin{exmp}\label{exmp.1}  Let $X$ be a Fano 3-fold  with an action  $\PGL_2 \curvearrowright X$ with a dense orbit. The complement is an ample divisor $D$.
    Thus  $\aut^\circ\bigl(\p_X\bigl(\o_X{\oplus}\o_X(nD)\bigr)\bigr)$ is not contained in any maximal, connected,   algebraic subgroup of
  $\bir(\p^4)$ for $n\geq 1$.

  There are only 2 such examples, both described in \cite{MR0726439}. One is now called the  {\it Mukai-Umemura} 3-fold  $X \subset \p^{21}$. 
  The other  is the unique $X\subset \p^6$ whose general hyperplane section is a degree 5 Del~Pezzo surface.
\end{exmp}

\begin{exmp}\label{exmp.2}  Let $g(u_0,u_1)$ be a homogeneous form without multiple roots and 
  $$
  Q_g:=(x_0^2 - x_1x_2 - g(u_0,u_1)x_3^2 =0)\subset \p^3_x\times \p^1_u.
  $$
  Together with a certain  $\PGL_2 $-action on the $\p^3_x$-factor as in  \cite[4.4.4]{MR4537333}, these are called  {\it  Umemura quadric fibrations.}

  We again get that
  $\aut^\circ\bigl(\p_{Q_g}(\o_{Q_g}{\oplus} L)\bigr)$ is not contained in any maximal, connected,   algebraic subgroup of
  $\bir(\p^4)$.
  The higher dimensional Umemura  fibrations studied in \cite{floris2024umemura}  give
similar examples in $\bir(\p^n)$ for $n>4$. \end{exmp}

The subtle example of \cite{ffz}  can also be treated using Corollary~\ref{main.cor.1}.

\begin{exmp}\label{exmp.3}   Let $X$ be a non-rational  but rationally connected 3-fold.
  Then  $\{1\}$ is the maximal, connected,   algebraic subgroup of $\bir(X)$.  We get that $\aut^\circ\bigl(\p_X(\o_X{\oplus}L)\bigr)$ is not contained in any maximal, connected,   algebraic subgroup of
  $\bir(X\times \p^1)$.

  As noted in \cite{ffz}, there is such an $X$ for which
  $X\times \p^2\simb \p^5$; see \cite{MR86m:14009}.
  \cite{ffz} uses this to obtain subgroups
  $$
  G\cong \aut^\circ\bigl(\p_X(\o_X{\oplus}L)\bigr)\times \PGL_{n-3}\subset \bir(\p^n)
  $$ that are  not contained in any maximal  algebraic subgroup of
  $\bir(\p^n)$ for $n\geq 5$.

  The construction in  \cite{ffz} uses  a non-ample line bundle; these can be treated using  Theorem~\ref{main.G.thm}.

\end{exmp}

The next example shows that the assumption
$(X, L)\not\cong \bigl(\p^n, \o_{\p^n}(1)\bigr)$ is necessary in Theorem~\ref{main.thm}.
(Strictly speaking, $\o_{\p^n}(1)$ does not have a $\PGL_{n+1}$-linearization for $n>0$.
However,  is has an $\SL_{n+1}$-linearization, which is enough for 
Theorem~\ref{main.thm}.)

    \begin{exmp}\label{exmp.4}  Let 
          $X=\p^n$ and $L=\o_{\p^n}(1)$. Then
          $\p_{\p^n}\bigl(\o_{\p^n}{\oplus} \o_{\p^n}(1)\bigr)$ is birational to $\p^{n+1}$, hence
  $$
  \aut^\circ\bigl(\p_{\p^n}\bigl(\o_{\p^n}{\oplus} \o_{\p^n}(1)\bigr)\bigr)\subset
  \aut^\circ(\p^{n+1})\cong \PGL_{n+2}.
  $$
  By contrast, $\aut^\circ\bigl(\p_{\p^n}\bigl(\o_{\p^n}{\oplus} \o_{\p^n}(r)\bigr)\bigr)$
  is a maximal, connected subgroup of $\bir(\p^{n+1})$ for $r\geq 2$ by
  Corollary~\ref{main.cor.2}  (though this is also easy to see directly).  
        \end{exmp}

\section{Groups acting on ruled varieties}

The idea behind the proof could be summarized as follows.

\begin{principle} Let $G\subset \bir(Y)$ be  an algebraic subgroup.
    If a normal,  unipotent subgroup  $ U\subset G$ acts on $Y$ with small orbits, then
    every  $G\subset \bar G\subset \bir(Y)$ is likely to have
a normal,  unipotent subgroup  $ \bar U\subset \bar G$  that acts on $Y$ with the same   orbits.
 \end{principle}

Theorem~\ref{main.G.thm} illustrates this  when   the orbits are 1 dimensional.

     \begin{defn} Let
       $
       \bir\bigl(X\times (\p^1, \infty)\to X\bigr)\subset \bir(X\times \p^1)
       $
       denote the subgroup of those birational maps $\phi:X\times \p^1\map X\times \p^1$ that commute with projection to $X$ and  map $X\times \{\infty\}$ birationally to $X\times \{\infty\}$. Thus there is a birational map
       $\phi_X:X\map X$ and a commutative diagram
       $$
       \begin{array}{ccc}
         X\times \p^1 & \stackrel{\phi}{\map} & X\times \p^1\\
         \pi\downarrow\hphantom{\pi} && \hphantom{\pi}\downarrow\pi \\
         X & \stackrel{\phi_X}{\map} & X,
       \end{array}
       $$
       where $\pi$ is the first coordinate projection.

       Note that $\phi\mapsto \phi_X$ is a group homomorphism; its kernel is denoted by  $\bir_X\bigl(X\times (\p^1, \infty)\bigr)$. Thus we have a group extension
       $$
       1\to \bir_X\bigl(X\times (\p^1, \infty)\bigr)\to
       \bir\bigl(X\times (\p^1, \infty)\to X\bigr)\to
       \bir(X)\to 1.
       $$
       \end{defn}

     \begin{defn}
       Given any connected, algebraic, subgroup
       $G\subset \bir\bigl(X\times (\p^1, \infty)\to X\bigr)$,
       let $G^v:=G\cap \bir_X\bigl(X\times (\p^1, \infty)\bigr)$ denote its
       {\it vertical subgroup,} and $G^h:=G/G^v\subset \bir(X)$
       its {\it horizontal quotient.}
       Thus we have a group extension
       $$
       1\to G^v\to G\to G^h\to 1.
       $$
       Note that $G^v$ is a solvable, normal subgroup of $G$.
   Let $U^v_G$ denote its  unique maximal, connected, unipotent subgroup. It is a normal subgroup of $G$.
   Let $U_G\subset G$ denote a maximal, connected, unipotent subgroup; it is unique up to conjugation.
   Note that  $U^v_G\subset U_G$,  and $U^h_G:=U_G/U^v_G\subset G^h$ is a
   maximal, connected, unipotent subgroup.
   (Connectedness is automatic in characteristic 0.)
        \end{defn}

     The following is the main result about subgroups of
     $ \bir(X\times \p^1)$. 
       
     \begin{thm}\label{main.G.thm}
       Let $X$ be a normal, projective variety and
        $G\subset \bir\bigl(X\times (\p^1, \infty)\to X\bigr)$
        a connected, algebraic subgroup. Assume that
       \begin{enumerate}
      \item $G^h\subset \bir(X)$ is a maximal, connected, algebraic subgroup,
      \item  $\dim U^v_G\geq 2$, and
      \item  for every  connected, algebraic subgroup  $ \{1\} \neq W\subset U^h_G$, the fixed point set of the  $W$-action by conjugation on $U^v_G$ has
    codimension $\geq 2$ in  $U^v_G$.
\end{enumerate}
       Then every connected, algebraic subgroup
       $$
       G\subset \bar G\subset \bir(X\times \p^1)
       $$
       is contained in $\bir\bigl(X\times (\p^1, \infty)\to X\bigr)$.
      Furthermore
        \begin{enumerate}\setcounter{enumi}{3}
        \item  $\bar G^h=G^h$, so $U^h_{\bar G}=U^h_G$, 
        \item $G^v\subset \bar G^v$, and
        \item $U^v_G\subset U^v_{\bar G}$.
        \end{enumerate}
        In particular, $\bar G$ also satisfies (\ref{main.G.thm}.1--3).
        \end{thm}

     {\it Note.}  (\ref{main.G.thm}.3) implies (\ref{main.G.thm}.2) whenever
     $ \{1\} \neq U^h_G$.
     \medskip

     Proof. We may assume that  $G\neq \bar G$.

     First we show that there is an intermediate group
     $
          G \subsetneq G_1 \subset \bar G,
          $
          such that
          $G_1\subset \bir\bigl(X\times (\p^1, \infty)\to X\bigr)$.
         
          To see this,  we apply the   purely group theoretic Lemma~\ref{group.00.lem}
           with $R:=U^v_G$.
           Assumption (\ref{main.G.thm}.3) says 
           that the alternative
           (\ref{group.00.lem}.2) does not happen.

          Thus there is a subgroup
          $\{1\}\neq R'\subset U^v_G$ whose (connected) normalizer
          $N_{\bar G}(R')^\circ$ is not contained in $G$.
          
          A key point is that generically $ R'$ and $U^v_G$ have the same orbits. In particular, 
          the closures of the 1-dimensional $R'$-orbits are still the
          fibers of the projection $\pi:X\times \p^1\to X$. Thus
          $N_{\bar G}(R')^\circ$ birationally preserves $\pi:X\times \p^1\to X$.
          The section $X\times \{\infty\}$ is the fixed point set of the $R'$-action, so this is also preserved. Thus
          $$
          N_{\bar G}(R')^\circ\subset \bir\bigl(X\times (\p^1, \infty)\to X\bigr).
          $$
           We can then take
           $
           G_1:=\langle G, N_{\bar G}(R')^\circ\rangle\subset \bar G$.

     Now $G_1$ clearly satisfies  $G^v\subset G_1^v$, 
         $U^v_G\subset U^v_{G_1}$, and
     $G^h\subset G_1^h\subset \bir(X)$. Since
     $G^h$ is a maximal, connected, algebraic subgroup, we must have $G^h= G_1^h$.
      Therefore $G_1$ also satisfies (\ref{main.G.thm}.1--3), and 
           induction on $\dim \bar G-\dim G$  completes the proof.
           \qed

           \begin{rem} If we do not assume that $G^h$ is maximal, the above proof still gives a $G\subsetneq G_1\subset \bar G$  that is a subgroup of
             $\bir\bigl(X\times (\p^1, \infty)\to X\bigr)$.
             However, if $U_G\neq U_{G_1}$ then 
             (\ref{main.G.thm}.3) may not hold for $G_1$, and the induction may break down.
\end{rem}

\section{Lemmas on algebraic groups}

\begin{notation}  For an algebraic group $G$, we let $G^\circ$ denote the maximal, connected subgroup.
              $B_G\subset G$ denotes a Borel  subgroup, 
         $U_G\subset B_G$ a maximal, connected,  unipotent subgroup,
  and  $R_G\subset U_G$ the unipotent radical of $G$.
  $Z(\ )$ denotes the center,  and $N_G(\ )$ the normalizer in $G$.
\end{notation}

\begin{lem}\label{group.0.lem} Let $G\subsetneq G'$ be connected, algebraic groups.
Then there is a subgroup $V\subset U_G$ of codimension $\leq 1$ such that
 $N_{G'}\bigl(Z(V)\bigr)^\circ\not\subset G$.
\end{lem}

Proof. Let $U_G\subset U_{G'}$  be maximal unipotent subgroups.

If $U_G\neq U_{G'}$ then let $U_G\subset N\subset U_{G'}$ be the normalizer of $U_G$ in $U_{G'}$.  Note that $U_G\neq N$ since $U_{G'}$ is nilpotent, and
$N\not\subset G$ since $U_G$ is a  maximal unipotent subgroup of $G$.
Then $N\subset N_{G'}(Z(U_G))^\circ$, and we can take $V:=U_G$.

Otherwise  $U_G= U_{G'}$.
For every simple root of $G'$ there is a minimal parabolic subgroup
$B_{G'}\subsetneq P_i\subset G'$; let
$V_i\subset U_G$ be its unipotent radical.
Then $P_i/V_i$ is either  $\SL_2$ or $\PSL_2$.
Therefore  $V_i\subset U_{G'}= U_G$ has codimension 1 and
$P_i\subset N_{G'}(V_i)$.

The $P_i$  generate $G'$, so  $P_i\not\subset  G$ for some $i_0$.
Then $V:=V_{i_0}$ works.  \qed

\begin{lem}\label{group.00.lem} Let $G\subsetneq G'$ be  connected, algebraic groups,
  $U_G\subset G$ a maximal, connected,  unipotent subgroup
  and  $\{1\}\neq R\subset U_G$ a connected, normal subgroup.
  Then one of the following holds.
  \begin{enumerate}
  \item  There is a connected subgroup $\{1\}\neq R'\subset R$  such that
 $N_{G'}(R')^\circ\not\subset G$.
  \item  There is a connected subgroup  $\{1\}\neq W\subset U_G/R$ such that 
the  $W$-action by conjugation on $R$ is 
    trivial on a subgroup $V_R\subset R$ of codimension $\leq 1$ in $R$.
   \end{enumerate}
\end{lem}

Proof. By (\ref{group.0.lem})  there is a subgroup $V\subset U_G$ of codimension $\leq 1$ such that
$N_{G'}\bigl(Z(V)\bigr)^\circ\not\subset G$. We are done if
$Z(V)\subset R$.
Otherwise  $$W:=\langle Z(V), R\rangle/R\cong Z(V)/\bigl( Z(V)\cap R\bigr)$$ acts trivially    on $V_R:=V\cap R$, which has codimension $\leq 1$. \qed

\section{Regularizing group actions}

First recall some  results on automorphisms of vector bundles; see for example \cite{MR0280493} or \cite[Sec.3]{ffz}. 
    
  \begin{say}[Automorphisms of vector bundles]\label{aut.vb.say}
    Let $C$ be a smooth, projective curve and $E$ a rank 2 vector bundle over $C$. Up to twisting with a line bundle, we have one of the following cases.
\begin{enumerate}
\item   $E$ is stable,  and then $\aut_C(E)\cong \gm$.
\item $E\cong \o_C{\oplus}\o_C$, and then $\aut_C(E)\cong \GL_2$.
  \item There is  an  extension $0\to L \to E\to \o_C\to 0$,
  where $\deg L\geq 0$ and the extension is non-split if $L\cong \o_C$.
  Then the unipotent radical of $\aut(E)$ is naturally isomorphic to $ H^0(C, L)$, and
  \begin{enumerate}
\item  $\aut_C(E)/H^0(C, L)\cong \gm$ if the extension is non-split, and
\item  $\aut_C(E)/H^0(C, L)\cong \gm^2$ if the extension is split.
\end{enumerate}
\end{enumerate}
Let now $\bar X$ be a normal, projective variety and
$X\subset \bar X$ an open subset such that $\bar X\setminus X$ has codimension $\geq 2$. Let $E$ be a rank 2 vector bundle over $X$ such that  $\aut(E)$ contains a unipotent subgroup of dimension $\geq 2$.
Working on general curve sections and using (\ref{aut.vb.say}.1--3) we get that, after twisting with a line bundle and possibly shrinking $X$, there is  a unique  extension
$$
0\to L \to E\to \o_X\to 0
\eqno{(\ref{aut.vb.say}.4)}
$$
such that $h^0(X, L)\geq 2$.
Furthermore
 \begin{enumerate}\setcounter{enumi}{4}
   \item the extension (\ref{aut.vb.say}.4) is $\aut_X(E)$ invariant, 
   \item the unipotent radical of $\aut_X(E)$ is naturally isomorphic to $ H^0(X, L)$, 
\item  $\aut_X(E)/H^0(X, L)\cong \gm$ if (\ref{aut.vb.say}.4) is non-split, and
\item  $\aut_X(E)/H^0(X, L)\cong \gm^2$ if (\ref{aut.vb.say}.4) is split.
\end{enumerate}

  \end{say}

  \begin{say}[Automorphisms of unstable $\p^1$-bundles]\label{aut.pb.say}
    Continuing with the notation of (\ref{aut.vb.say}), 
     let $\pi:P:=\p_X(E)\to X$ be the corresponding $\p^1$-bundle.
     Then $\aut_X(P)$ is the quotient if $\aut_X(E)$ by the scalars.
    Thus we get that, if 
 $\aut_X(P)$ contains a unipotent subgroup of dimension $\geq 2$, then 
 \begin{enumerate}
   \item there is a unique $\aut_X(P)$ invariant  section $D\subset P$, 
   \item the unipotent radical of $\aut_X(P)$ is naturally isomorphic to $ H^0(X, L)$, 
\item  $\aut_X(P)=H^0(X, L)$ if (\ref{aut.vb.say}.4) is non-split, 
\item  $\aut_X(P)/H^0(X, L)\cong \gm$ if (\ref{aut.vb.say}.4) is split, and
\item  the splittings of $E\to \o_X$ correspond to the
  splittings of  $\aut_X(P)\to \gm$. 
\end{enumerate}
 Note that the natural way to get a vector bundle on $X$  is
 to take  $\pi_*\o_P(D)$. This sits in  the exact sequence
 $$
 0\to \o_X=\pi_*\o_P\to \pi_*\o_P(D) \to \o_P(D)|_D\to 0.
 \eqno{(\ref{aut.pb.say}.6)}
 $$
 Thus $\pi_*\o_P(D)\cong E\otimes L^{-1}$.

    \end{say}
 
Next we show  a converse to these observations.

\begin{prop} \label{best.reg.G.prop}
  Let $\bar X$ be a normal, projective variety with smooth locus $X\subset \bar X$, and $G\subset \bir\bigl(X\times (\p^1, \infty)\to X\bigr)$ an algebraic subgroup. Assume that $G^h\subset \aut(\bar X)=\aut(X)$. Then
there is a $G^h$-invariant open subset   $X^\circ\subset X$,
  and a $\p^1$-bundle $\pi^\circ: P^\circ\to X^\circ$ with a rational section
  $D^\circ\subset P^\circ$, such that
  $$
  G\subset \aut\bigl((P^\circ, D^\circ)\to X^\circ\bigr).
  $$
  Moreover
  $X\setminus X^\circ$  can be chosen to have  codimension $\geq 3$ in $X$.
\end{prop}

\begin{thm} \label{best.reg.G.prop.2} Using the notation and assumptions of
  (\ref{best.reg.G.prop}), assume in addition that there is a $\gm\cong T\subset G^v/U^v_G$  and $\dim U^v_G\geq 1$. 
  Then there is a unique $\p^1$-bundle 
  $\pi^m: P^m\to X$ with a  $G$-invariant section $D^m\subset P^m$,  such that
    \begin{enumerate}
     \item the
      set of fibers with trivial $U^v_G$-action has codimension $\geq 2$, and
    \item for every other $\pi': (P', D') \to X'$
satisfying $G\subset \aut\bigl((P', D')\to X'\bigr)$, 
      there is a natural $G$-equivariant 
      $$
      \pi'_*\o_{P'}(D')\into \pi^m_*\o_{P^m}(D^m).
      $$
      \end{enumerate}
  \end{thm}

 Proof of (\ref{best.reg.G.prop}) and (\ref{best.reg.G.prop.2}).  By \cite{MR0074083} the $G$-action  is regularized on some $Y_1\simb \bar X\times \p^1$.
    Let $Y_2\subset \bar X\times Y_1$ be the closure of the graph of
    the rational map $Y_1\map \bar X\times \p^1\to \bar X$.
    Then $G$ acts regularly on $Y_2$.  After normalization we get a  $G$-equivariant  $Y_3\to \bar X$, where $Y_3$ is normal.
    It is  a $\p^1$-bundle over a $G^h$-invariant, dense open subset
    $X_3\subset X$.

    Let $x_i\in X\setminus X_3$ be the  generic points that have codimension 1 in $X$.
    Working over the localization $\o_{x_i,X}$, (\ref{over.dvr.1}) gives a
    $G$-invariant $\p^1$-bundle extending $Y_3\to X_3$ over  some open neighborhood of $x_i$.
    Gluing these to $Y_3\to X_3$ we get a $G$-equivariant $\p^1$-bundle
    $\pi_4: Y_4\to X_4$ with a rational section $D_4\subset Y_4$,
    where   $X\setminus X_4$  has codimension $\geq 2$ in $X$.

    Let $F$ denote the  push-forward of $(\pi_4)_*\o_{Y_4}(D_4)$  to $X$.
    Let $X^\circ\subset X$ be the largest open subset where $F$ is locally free. Then $X\setminus X^\circ$  has codimension $\geq 3$ in $X$.

    Set $P^\circ:=\p_{X^\circ}\bigl(F|_{X^\circ}\bigr)$.
    The $G$-action  extends to a regular $G$-action on
    $P^\circ$ by  by the Matsusaka-Mumford lemma; see
    \cite[11.39]{k-modbook} for the version that we need.
    This completes the proof of 
    (\ref{best.reg.G.prop}).

In  (\ref{best.reg.G.prop.2}) the  $T$-action on $(\pi_4)_*\o_{Y_4}(D_4)$ produces a splitting
$
(\pi_4)_*\o_{Y_4}(D_4)\cong \o_{X_4}\oplus   L^{-1}_4,
$ and
    $U^v_G$ is naturally identified with a subspace
$W_G\subset H^0(X, L_4)$.
Let $L_5\subset L_4$ be the invertible subsheaf generated by the sections in
$W_G$ at codimension 1 points.  Then $L_5$ extends to a line bundle $L$ on $X$,
and $G$ also acts on
$E:=\o_{X}\oplus   L^{-1}$. Its projectivization 
gives the reguired
$\p^1$-bundle for  (\ref{best.reg.G.prop.2}).

Uniqueness can be checked on curves. We can thus work with a rank 2 free module
$E_R:=Re_1\oplus Re_2$ over a DVR $R$ with maximal ideal $(t)$ and quotient field $K$, where
the $T$ action is  $(e_1, e_2)\mapsto (\lambda e_1, \lambda^{-1}e_2)$,
and the elements of the unipotent subgroup act as  
$$
(r_1e_1, r_2e_2)\mapsto \bigl(r_1e_1, (r_2+ur_1) e_2\bigr) \qtq{where} u\in R.
$$
The  unipotent action  is nontrivial on the central fiber iff not all the $u$ are divisible by $t$.

Every other $T$-equivariant extension of  $E_K$ is of the form
$E^*_R:=t^{c_1}Re_1\oplus t^{c_2}Re_2$. In the basis $f_i:=t^{c_i}e_i$ the  unipotent action   becomes
$$
\bigl(r_1f_1, r_2f_2\bigr)\mapsto
\bigl(r_1e_1,  (r_2+ut^{c_1-c_2}r_1) f_2\bigr).
$$
This is a regular action iff $c_1\geq c_2$ and  it is
nontrivial on the central fiber iff $c_1=c_2$.
In this case the projectivizations of  $E_R$ and of $E^*_R$ are naturally isomorphic. \qed

    \begin{rem} If $\bar X$ is factorial, then $P^m$ in
      (\ref{best.reg.G.prop.2}) extends to a $\p^1$-bundle.
      If $\bar X$ is $\q$-factorial, then  $L$ extends  to a $\q$-line bundle $\bar L$, and 
      $P^m$ extends to an orbifold  $\p^1$-bundle 
$$
  \proj_{\bar X}\Bigl(\oplus_{r\geq 0}\bigl(\oplus_{i=0}^r \bar L^{[i]}\bigr)\Bigr).
      $$
      \end{rem}

    \begin{exmp} The codimension $3$ subset $X\setminus X^\circ$ is usually unavoidable in  (\ref{best.reg.G.prop}), as shown by the following example taken from \cite{MR1620110}.  Consider  the second coordinate projection  
    $$
    \pi: Y:=(\tsum_i x_iy_i=0)\subset \p^2_{y}\times \a^3_{x}\to \a^3_{x}.
    $$
    The fiber over $(0,0,0)$ is $\p^2$, all other fibers are $\p^1$.
    We can chose $D:=(y_0=0)$ as a rational section.
    Note that $Y$ is smooth and 
the projection $\pi$ is $\GL_3$-equivariant for the action
    $
    \bigl({\mathbf x}, {\mathbf y}\bigr)
    \mapsto
    \bigl(A^t{\mathbf x}, A^{-1}{\mathbf y}\bigr).
    $
    (This shows that \cite[Lem.2.2]{ffz} is not quite true.)
    \end{exmp}

     \begin{lem}\label{over.dvr.1}  Let $k$ be a field, $(t\in T)$ the spectrum of an excellent DVR over $k$, and  $\pi:P\to T$ a projective morphism with general fiber $\p^1$, and
       $D\subset P$ a section.  Let $B$ be the group of $k$-automorphisms
       of $\pi:(P,D)\to T$.

       Then $\pi:(P,D)\to T$ is $B$-equivariantly birational to a trivial
       bundle  $\pi': (\p^1\times T, \{\infty\}\times T)\to T$.
     \end{lem}

     Proof. We can resolve the singularities by normalizing and blowing up the singular set, hence the $B$ action the lifts. Thus we may assume that $P$ is regular.
     Note that  $(-1)$-curves in the central fiber are disjoint from each other, we can thus contract all $(-1)$-curves disjoint from $D$ in a  $B$-equivariant way.
     After such contractions we get a minimal conic bundle with a section, hence a $\p^1$-bundle.  (The residue field of $T$ may be imperfect, see \cite{MR3795482} for a proof in these cases.) \qed

     \begin{exmp} The non-uniqueness without the $\gm$-action is illustrated by the following, at the vector bundle level.

       Let $R$ be a DVR with maximal ideal $(t)$.
       Consider the rank 2 free module $E:=Re_1{\oplus}Re_2$ with an endomorphism
       given by $\phi(e_1)=e_2, \phi(e_2)=0$. The $\phi$-action is nontrivial on $E/tE$. Another extension of $Re_1$ by $Re_2$ is given by 
       $E':=R(e_1+t^{-1}e_2){\oplus}Re_2$. The $\phi$-action is again nontrivial on $E'/tE'$.
       \end{exmp}

     \section{Proof of Theorem~\ref{main.thm}}

     \begin{say}[Proof of Theorem~\ref{main.thm}]\label{main.thm.pf}
       First we aim to apply (\ref{main.G.thm}).
       Note that (\ref{main.G.thm}.1) is assumed in Theorem~\ref{main.thm}.
       There is nothing to prove if $X$ is a point, otherwise (\ref{main.G.thm}.2) also holds. 
     We check in (\ref{exception.lem}) that     (\ref{main.G.thm}.3) 
     is satisfied.
     Thus we get that
     $$
      G\subset \bir\bigl(X\times (\p^1, \infty)\to X\bigr).
     $$
     Now we can apply (\ref{best.reg.G.prop.2}) to get that
     $$
           G\subset \aut^\circ\bigl(\p_X\bigl(\o_X{\oplus}L_G\bigr)\bigr)
           $$
           for some line bundle $L_G$.

           Set $G_0:=\aut^\circ\bigl(\p_X\bigl(\o_X{\oplus}L\bigr)\bigr)$.
           Since $L$ is generated by its global sections,
           $\p_X(\o_X{\oplus}L)\to X$ is isomorphic to the regularization of $G_0$ that is denoted by  $P^m$ in (\ref{best.reg.G.prop.2}.2).
The inclusion in (\ref{best.reg.G.prop.2}.2) induces  a  $G$-equivariant  $L_G^{-1}\into L^{-1}$, thus $L_G=L(D)$ for some effective  $G^h$-invariant divisor $D$. \qed  
           
          \end{say}
     
The following accounts for the exceptional case  $(X, L)\not\cong \bigl(\p^n, \o_{\p^n}(1)\bigr)$
in  Theorem~\ref{main.thm}.

    \begin{lem} \label{exception.lem} Let $X$ be a normal, projective variety with a nontrivial $\ga$-action.
      Let $L$ be a very ample, $\ga$-linearized line bundle.
      Assue that the induced $\ga$-action is trivial on a codimension 1 subspace  $W\subset H^0(X, L)$.

Then $X$ is a cone, and the  $\ga$-action is along its  generators.

If $X$ is smooth of dimension $n$,  then $(X,L)\cong \bigl(\p^n, \o_{\p^n}(1)\bigr)$.
    \end{lem}

    Proof.
    Let $C\subset X$ be the closure of a $\ga$-orbit and
    $c\in C$ the fixed point.
$L$ gives an imbedding $X\into \p^N$; let $\langle C\rangle\subset \p^N$ be the linear span of $C$. 
    Let $\pi:\bar C\to C$ be the normalization.
    Then
    $$
    \dim \im\bigl[H^0(X, L)\to H^0(\bar C, \pi^*L)\bigr]= 1+\dim \langle C\rangle.
    $$
    The $\ga$-action  on $H^0\bigl(\p^1, \o_{\p^1}(d)\bigr)$ is trivial only on a 
    1-dimensional  subspace.  Thus the image of $W$ in
    $H^0(\bar C, \pi^*L)$ has dimension 1 and codimension 1.
    Therefore $\dim \langle C\rangle=1$, thus $C$ is a line,  and the  image of $W$ consists of the sections that vanish at $c$.

Given 2 orbits  $C_i\ni c_i$, if $c_1\neq c_2$ then vanishing at both gives a subset of codimension $\geq 2$. Thus all $\ga$-orbit  closures are lines, and 
all these lines pass trough the same point. So $X$ is a cone.

The only smooth cone is a linear subspace.
\qed

\begin{ack} The main impetus came from the lecture {\it On Algebraic Subgroups of the Cremona Group} of Floris at the  Conference on Higher Dimensional Geometry (May, 2024) at the
  Simons Center for   Mathematics \& Physical Sciences.
  
  In the lecture Floris explained  the main steps very clearly, but omitted the technical details. So, when I was trying to think of a proof after the talk, I was free to start on  a  somewhat different path.

  I thank E.~Floris and P.~Fong for helpful comments.
  Partial  financial support    was provided  by  the NSF under grant number
DMS-1901855.
\end{ack}


\def\cprime{$'$} \def\cprime{$'$} \def\cprime{$'$} \def\cprime{$'$}
  \def\cprime{$'$} \def\dbar{\leavevmode\hbox to 0pt{\hskip.2ex
  \accent"16\hss}d} \def\cprime{$'$} \def\cprime{$'$}
  \def\polhk#1{\setbox0=\hbox{#1}{\ooalign{\hidewidth
  \lower1.5ex\hbox{`}\hidewidth\crcr\unhbox0}}} \def\cprime{$'$}
  \def\cprime{$'$} \def\cprime{$'$} \def\cprime{$'$}
  \def\polhk#1{\setbox0=\hbox{#1}{\ooalign{\hidewidth
  \lower1.5ex\hbox{`}\hidewidth\crcr\unhbox0}}} \def\cdprime{$''$}
  \def\cprime{$'$} \def\cprime{$'$} \def\cprime{$'$} \def\cprime{$'$}
\providecommand{\bysame}{\leavevmode\hbox to3em{\hrulefill}\thinspace}
\providecommand{\MR}{\relax\ifhmode\unskip\space\fi MR }
\providecommand{\MRhref}[2]{%
  \href{http://www.ams.org/mathscinet-getitem?mr=#1}{#2}
}
\providecommand{\href}[2]{#2}

  \bigskip

  Princeton University, Princeton NJ 08544-1000, \

  \email{kollar@math.princeton.edu}

     \end{document}